\documentclass{amsart}
\font\nsc=cmcsc10 scaled 900

\let\<=\langle
\let\>=\rangle
\let\noi=\noindent
\let\sse=\subseteq
\let\limply=\Longrightarrow

\def\0{\{0\}}
\def\smallfrac#1#2{{\textstyle{\frac{#1}{#2}}}}
\def\esup{{\kern1pt{\rm ess}\kern.5pt\sup\kern1.5pt}}
\def\notlimply{{\,{{\limply}\kern-11pt{\slash}}\;\;\;}}
\def\noimply{{\kern3pt\not\kern-3pt\limply}}

\def\C{{\kern1pt\mathcal C}}
\def\H{{\mathcal H}}
\def\M{{\mathcal M}}
\def\N{{\mathcal N}}
\def\R{{\mathcal R}}
\def\X{{\mathcal X}}

\def\CC{{\mathbb C\kern.5pt}}
\def\NN{{\mathbb N\kern.5pt}}
\def\TT{{\mathbb T\kern.5pt}}

\newsymbol\blacksquare 1004

\def\smallmatrix#1{\null\,\vcenter{
                    \baselineskip=8pt\mathsurround=0pt\ialign{
                    \hfil ${\scriptstyle##}$
                    \hfil &&
                    \hfil ${\scriptstyle##}$
                    \hfil \crcr
                    \mathstrut \crcr
                    \noalign{\kern-\baselineskip}#1 \crcr
                    \mathstrut \crcr
                    \noalign{\kern-\baselineskip} \crcr }}\!}

\begin{document}

\vglue-70pt\noi
\hfill{\it }\phantom{{\bf XX} (20XX) xx--xx (XXXX) to appear}

\vglue20pt
\title[Normaloid Operators and the Root Problem]
      {Normaloid Operators and the Root Problem}
\author{B.P. Duggal}
\address{Faculty of Sciences and Mathematics, University of Ni\v s, Ni\v s,
Serbia}
\email{bpduggal@yahoo.co.uk}
\author{C.S. Kubrusly}
\address{Catholic University of Rio de Janeiro, Rio de Janeiro, Brasil}
\email{carlos@ele.puc-rio.br}
\author{H.M. Stankovi\' c}
\address{Faculty of Electronic Engineering, University of Ni\v s, Ni\v s,
Serbia}
\email{hranislav.stankovic@elfak.ni.ac.rs}
\renewcommand{\keywordsname}{Keywords}
\date{June 14, 2026}

\vskip-10pt
\begin{abstract}
The paper extends previous results on the $n$th root problem to a large
class of Hilbert-space operators, namely, the class of all normaloid operators
with normaloid parts, which includes the paranormal operators, and also the
$k$-paranormal operators.\ It is shown that if a normaloid operator with
normaloid parts has a normal $n$th power, then it is normal.\
\end{abstract}

\maketitle

\vskip-10pt\noi
\section{Introduction}

The $n$th root problem asks what classes of Hilbert-space operators
$T\kern-1pt$ are such that, if $T^n\!$ is~normal, then $T\kern-1pt$ is
normal.\ The problem has been investigated for more than six decades and
remains under active research these days.\ Perhaps the most popular example
along this line is the celebrated result for hyponormal operators proved by
Stampfli in the early sixties \cite[Theorem 5]{Stampfli}, which was extended
to paranormal operators by Ando a decade afterwards \cite[Theorem 6]{And}.\
The problem has recently been extended to $k$-paranormal operators and beyond
\cite[Theorem 3.1]{SK1} and \cite[Theorem 4.1]{SK2} (and also the references
therein).\

\vskip6pt
The main results proved here are Theorems 5.2 and 6.1, which read as follows.\
\vskip4pt
\begin{description}
\item{$\kern-2pt\circ\kern2pt$}
If a Hilbert-space operator $\kern-.5pt T\kern-1pt$ is normaloid with an
invertible $n$th power~$\kern-.5pt T^n\kern-1pt$~such that
${\|T^n\kern-.5pt\|^{-1}\kern-2.5pt=\kern-1.5pt\|T^{-n}\kern-.5pt\|}$,
then $\kern-.5ptT\kern-1pt$ is a multiple of a unitary
operator~(\hbox{Theorem}~5.2).\
\vskip4pt
\item{$\kern-2pt\circ\kern2pt$}
If an operator $T\kern-1pt$ on a Hilbert space is normaloid with normaloid
direct summands, and if $T^n\kern-1pt$ is normal for some positive integer
$n$, then $T\kern-1pt$ is normal (Theorem 6.1).\
\end{description}
\vskip-3pt

\vskip6pt
Theorem 5.2 will contribute to the proof of Theorem 6.1.\ The class of
normaloid Hilbert-space operators with normaloid direct summands is quite
large, being near to the class of normaloid operators themselves.\ It includes
the $k$-paranormal operators (and so, in particular, the classes of paranormal
and hyponormal operators).\

\vskip6pt
The paper is split into 7 sections.\ Basic terminology and notation are posed
in Section 2.\ The class of completely normaloid operators is introduced in
Section 3, and the root problem is considered in Section 4.\ Section 5
discusses in detail a class of normaloid operators with a norm condition on
their powers that become multiples of unitary operators.\ The root problem for
completely normaloid operators~is~inves\-tigated in Section 6.\ Additional
results close the paper in Section~7.\

\section{Basic Terminology and Notation}

By an operator we mean a bounded linear transformation of a normed space~$\X$
into itself.\ A subspace $\M$ of $\X$ is a closed linear manifold of $\X$;
it~is nontrivial if ${\0\ne\M\ne\X}.$ A subspace $\M$~is invariant for an
operator $T\kern-1pt$ (or is $T\kern-1pt$-invariant) if ${T(\M)\sse\M}.$ The
induced uniform norm of an operator $T\kern-1pt$ is denoted by $\|T\|$, and,
if $\X$ is a (complex) Banach space, the spectrum of $T\kern-1pt$ is denoted
by $\sigma(T)$ and its spectral radius by $r(T).$ If $T\kern-1pt$ is an
operator on a Hilbert space $\H$, then~the~operator $T^*\kern-1pt$ on $\H$
stands for the adjoint of $T\kern-1pt.$ A subspace $\M$ of an operator
$T\kern-1pt$ on a Hilbert space $\H$ reduces $T\kern-1pt$ (or $\M$ is a
reducing subspace for $T\kern-1pt$) if $\M$ is invariant for $T\kern-1pt$ and
for $T^*\!.$ Equivalently, if $\M$ and $\M^\perp\!$ are both
$T\kern-1pt$-invariant, where $\M^\perp\!$ stands for the orthogonal
complement of $\M$ in $\H.$ The restriction $T|_\M\!$ of $T\kern-1pt$ to a
reducing~subspace $\M$ is referred to as a direct summand of $T\kern-1pt$ or
as a {\it part}\/ of $T\kern-1pt.$ (Sometimes~the~term ``part'' is defined as
the restriction of an operator to an invariant subspace, but here we use the
term as a synonym of direct summand.) A part $T|_\M$ of an \hbox{operator}
$T\kern-1pt$~is nontrivial if $\M$ is a nontrivial reducing subspace for
$T\kern-1pt$, and $T\kern-1pt$ is reducible~if~it~has~a nontrivial part;
otherwise it is irreducible.\ The symbol $\oplus$ stands for the
\hbox{orthogonal} direct sum of subspaces, as in
${\H=\kern-1pt\M\oplus\M^\perp\!}$, or of operators, as in
${T\kern-1pt=T|_\M\oplus T|_{\M^\perp\!}}\!$ if $\M$ reduces $T\kern-1pt$,
where $T|_\M$ and $T|_{\M^\perp}\!$ are direct summands, that is, parts, of
$T\kern-1pt.$

\vskip6pt
A nonzero operator $T\kern-1pt$ on a linear space is nilpotent of index $j$ if
${T^j\kern-1pt=O}$ for some integer ${j>1}$ and ${T^i\!\ne O}$ for every
positive integer ${i<j}$, where $O$ denotes the null operator.\ A nilpotent
operator without a specified index will be supposed to be of index 2.\ An
involution is an invertible operator on a linear space for which
${T^{-1}\!=T}$; equivalently, for which ${T^2\!=I}$, where $I$ stands for the
identity operator.\

\vskip6pt
An operator $T\kern-1pt$ acting on a Hilbert space is normal if it commutes
with $T^*\kern-3pt$ (i.e., $T^*T\kern-1pt=T\kern1pt T^*\kern-1pt).$ It is
quasinormal if it commutes with ${T^*T\kern-1pt}$ (i.e., if
${(T^*T\kern-1pt-T\kern1pt T^*)\kern1pt T\!=O}).$ So normal operators are
quasinormal.\ A self-adjoint is a normal operator such that ${T=T^*\!}$, and a
unitary is an invertible normal operator such that ${T^{-1}\!=T^*}\!.$ The
following assertions are pairwise equivalent.\
\begin{description}
\item{$\kern3pt$\rm(i)$\kern7pt$}
$\,T$ is self-adjoint and unitary\/.
\vskip3pt
\item{$\kern1pt$\rm(ii)$\kern6pt$}
$\,T$ is a unitary involution\/.
\vskip3pt
\item{\rm(iii)$\kern4pt$}
$\,T$ is a self-adjoint involution\/.
\vskip3pt
\item{$\kern1pt$\rm(iv)$\kern4pt$}
$\,T$ is a normal involution\/.
\vskip3pt
\item{$\kern3pt$\rm(v)$\kern4pt$}
$\,T$ is a quasinormal involution\/.
\end{description}
(For the equivalence among items (i) to (iv), see, e.g.,
\cite[Problem 5.43]{EOT}; (iv) and (v) are equivalent because an invertible
quasinormal is normal.)\ An operator satisfying any of the above equivalent
conditions is called a symmetry.\

\section{Completely Normaloid Operators}

\vskip6pt
A Hilbert-space operator $T\kern-1pt$ is hyponormal if
${T\kern1pt T^*\kern-2pt\le T^*T}\!.$ (Quasinormal operators are hyponormal.)\
An operator $T\kern-1pt$ acting on a normed space $\X$ is paranormal if
$\|Tx\|^2\!\le{\|T^2x\|\kern.5pt\|x\|}$ and $k$-paranormal if
$\|Tx\|^{k+1}\!\le{\|T^{k+1}x\|\kern.5pt\|x\|^k}\kern-1pt$ for some positive
integer $k$, for every ${x\in\X}.$ A Banach-space operator $T\kern-1pt$ is
normaloid if its spectral radius coincides with its norm (i.e.,
${r(T)=\|T\|}$), which is equivalent to
${\|T^k\|\kern-1pt=\kern-.5pt\|T\|^k}\!$ for every integer
${k\kern-1pt\ge\kern-1pt1}.$ These classes are related by proper inclusion:\
$$
\hbox{\nsc Normal $\subset$ Hyponormal $\subset$ Paranormal $\subset$
$k$-Paranormal $\subset$ Normaloid}.                            \eqno{(\dag)}
$$
In fact, paranormal operators are normaloid \cite[Theorem 1(i)]{ISY}
(Example 4.1 below shows that the inclusion is proper).\ It is known that the
inverse of an invertible~para\-normal operator is paranormal
\cite[Theorem 1(ii)]{ISY}, that integer powers of paranormal operators are
paranormal \cite[Theorem 1]{Fur0} and \cite[Theorem 2, p.$\,$103]{Fur}, and
that the restriction of a paranormal operator to an invariant subspace is
paranormal \cite[p.$\,$153]{II}.\ It is also known that $k$-paranormal
operators are normaloid and that restrictions of $k$-paranormal operators to
invariant subspaces are $k$-paranormal (see, e.g.,
\cite[Proposition~1(b)]{KD1}$\kern.5pt$).\ Properties of powers and inverses
of $k$-paranormal operators have been investigated in, for instance,
\cite[$\kern-1pt$Theorems~1~and~2]{KD1}.\

\vskip6pt
A Banach-space operator is hereditarily normaloid if every restriction of it
to an invariant subspace is normaloid.\ Let {\sc HN} denote the class of all
hereditarily~normal\-oid operators.\ The class {\sc HN} was investigated in
\cite{DD} (see also \cite{DDK}, \cite{Dug2}$\kern.5pt$).\ This is~a~large
class of operators.\ It lies properly between the $k$-paranormal and the
normaloid~operators (see, e.g., \cite[Proposition 1]{KD1}$\kern.5pt$).\ We
will be dealing with a still larger class.\

\vskip6pt\noi
{\it Definition.}\/
A Hilbert-space operator is {\it completely normaloid}\/ if every
restriction of it to a reducing subspace is normaloid (i.e., if every part of
it is normaloid).\

\vskip6pt
Let {\sc CN} stand for the class of all completely normaloid Hilbert-space
operators.\ This class is trivially included in the class of normaloid
operators and includes the class of all irreducible normaloid operators.\
Then we get a \hbox{refinement}~of~chain~$(\dag)$:
$$
\hbox{\nsc $k$-Paranormal $\subset$ HN $\subset$ CN $\subset$ Normaloid}.
                                                               \eqno{(\ddag)}
$$
To verify that ${\hbox{\sc HN}\kern-1.5pt\sse\kern-1.5pt\hbox{\sc CN}}$, take
a Hilbert space $\H$ and an operator $T\kern-1.5pt$ on $\H.$
Let~$\hbox{\it Lat}\,(\H)$ be the lattice of all subspaces of $\H.$ Let
$\hbox{\sc NLD}$ denote the collection of all normaloid operators on any
${\!\M\in\hbox{\it Lat}\,(\H)}.$ Let
${\hbox{\it Lat}\,(T)\sse\hbox{\it Lat}\,(\H)}$ be the lattice of all
$T$-invariant subspaces for $T\kern-1pt.$ Let $\hbox{\it Red}\,(T)$ be the
lattice of all reducible subspaces~of~$T\kern-1pt.$
$$
T\in\hbox{\sc HN}
\iff
T|_\M\in\hbox{\sc NLD}\;\hbox{ \rm for all }\;\M\in\hbox{\it Lat}\,(T)
$$
$$
\limply\quad
T|_\M\in\hbox{\sc NLD}\;\hbox{ \rm for all }\;\M\in\hbox{\it Red}\,(T)
\sse\hbox{\it Lat}\,(T)
\iff
T\in\hbox{\sc CN}.
$$

\vskip6pt\noi
{\bf Remark 3.1.}\
In \cite[Definition 2]{II}, a Hilbert-space operator was said to be invariant
normaloid if the restriction of it to every invariant subspace is normaloid.\
(See~also \cite[Definition 7.1.15]{Ist}$\kern.5pt$).\ The notion of a
normaloid operator given in \cite[Definition 1]{II} was in terms of the
numerical radius (i.e., ${w(T)=\|T\|}$) rather than in terms of the spectral
radius (i.e., ${r(T)=\|T\|}$).\ These definitions of normaloidness, however,
are equivalent (see, e.g., \cite[Proposition 6.27]{EOT}$\kern.5pt$).\
Moreover, it was shown in \cite[Theorem 1]{II} that if an invariant normaloid
(i.e., a hereditarily normaloid) \hbox{Hilbert-space}~oper\-ator is compact,
then it is normal.\ This can be extended to completely~normaloid.\
$$
\hbox{\it If\/ ${T\in\hbox{\sc CN}}$ is compact, then it is normal}\/.
$$
In fact, the proof in \cite[Theorem 1]{II} actually shows that if a
Hilbert-space operator is compact and all its direct summands are normaloid,
then it is normal.\ Therefore, all subclasses of normaloid operators
considered here, when acting on a finite-dimen\-sional space, are normal.\
But there are nonnormal normaloid operators on~finite-dimensional spaces:\ for
instance, ${T\kern-2pt=\kern-1pt A\oplus I}\kern-.5pt$ with 
${A\kern-1pt=\kern-1pt\big(\smallmatrix{0   & 1 \cr
                                        1/2 & 0 \cr}\big)}.$
                                                         \hfill$\blacksquare$

\section{The Root Problem}

Let $\C$ be a class of operators on a Hilbert space of dimension greater than
one containing normal operators.\ Recall that powers and parts of a normal
operator are again normal.\ Let $n$ be a positive integer.\ The {\it $n$th
root problem}\/ (or the {\it root problem}\/, for short) asks for classes
$\kern-.5pt\C\kern-1pt$ of operators satisfying the following condition.\
\vskip6pt
\centerline{\it If\/ $T\kern-1pt$ lies in\/ $\C$ and\/ $T^n\kern-1pt$ is
normal, then\/ $T\kern-1pt$
is normal\/.}
\vskip6pt\noi
Since $T\kern-1pt$ belongs to $\C$, and since it is expected to be normal,
the class $\C$~is~supposed to contain normal operators; otherwise, the root
problem is vacuous.\

\vskip6pt
The functional calculi for operators acting on a complex Hilbert space ensure
that a normal operator $N$ has a normal $n$th root.\ In fact, if
${\psi\!:\CC\to\CC}$ is the function assigning to each
${\zeta\kern-1pt\in\CC}$ its principal $n$th root, that is,
${\psi(\zeta)=|\zeta|^{\frac{1}{\scriptstyle n}}
\hbox{exp}
\big(i\,\hbox{Arg}(\zeta)\kern1pt\frac{1}{\scriptstyle n}\kern-1pt\big)}$,
then ${\psi(\zeta)^n\kern-1pt=\zeta}$ for every ${\zeta\kern-1pt\in\CC}$ and
$\psi(N)$ is an $n$th root of $N$, which is again a normal operator.\ For more
aspects of the characterisation of the $n$th root problem, see, e.g.,
\cite{Dug,Ker}.\ However, there are two crucial features here.\ First, there
is no uniqueness:\ a normal operator may have several normal square roots.\
(Trivial example:\ $I$ and $-I$ are normal square roots of the identity
operator.)\ More \hbox{important}, there are square roots for normal operators
that are not normal.\ (Trivial example:~a~nonzero nilpotent operator, which is
never normal, is a square root of the null operator.)\

\vskip6pt\noi
{\bf Example 4.1.}\
An example next to trivial emphasises some critical questions.\ Take
$T\kern-1pt={L\oplus I}$, where $L$ is a nonzero nilpotent contraction.\ The
contraction $T\kern-1pt$ is not normal because $L$ is not normal; $L$ is not
even normaloid --- the only quasinilpotent normaloid operator is the null
operator.\ But ${T^2\!=O\kern-.5pt\oplus I}$ is a nonzero normal (thus not
nilpotent) with ${r(T^2)\kern-1pt=\kern-1pt1}$, and so $T\kern-1pt$ is
normaloid with ${\|T\|\kern-1pt=\kern-1ptr(T)\kern-1pt=\kern-1pt1}.$ A defect
of this example is that, although $T\kern-1pt$ is a reducible normaloid
operator, it has a nonnormaloid part $L$, thus dismissing any chance
of $T$ being~normal.\                                    \hfill$\blacksquare$

\vskip6pt
For comparison, consider a well-known class satisfying the root problem:\
\vskip6pt\noi
\centerline{\it If\/ $T\kern-1pt$ is hyponormal and\/ $T^n\kern-1.5pt$ is
normal for some\/ $n$, then\/ $T$ is normal\/}
\vskip6pt\noi
\cite[Theorem 5]{Stampfli}, which was extended to paranormal operators in
\cite[Theorem 6]{And} and to $k$-paranormal operators in
\cite[Theorem 3.1]{SK1}.\ (The root problem has been extended to
classes of normaloid --- and nonnormaloid --- operators that properly include
the paranormal operators; see, e.g., \cite{SK1,SK2,KS}.\ For variants~of~the
problem, for instance, involving subnormal and quasinormal operators together,
as well as multivariable versions of it, see, e.g., \cite{Sta,Sta2} and the
references therein.)\

\vskip6pt
Since the classes of hyponormal and paranormal operators (which are normaloid)
satisfy the root problem, they do not present the defect described in
Example 4.1.\ If a hyponormal (paranormal) is reducible, then its parts are
hyponormal (paranormal), and so its parts are normaloid.\ The same property
also applies to $k$-para\-normal operators (which are normaloid):\ parts of
$k$-paranormal are $k$-paranormal, thus normaloid --- cf.\ chains $(\dag)$ and
$(\ddag).$ Recalling that {\sc CN} stands for the class of normaloid operators
with normaloid parts, consider the following question.\

\vskip6pt\noi
{\bf Question 4.2.}\
{\it Take\/ $T\kern-1pt$ in\/ {\sc CN}.\ Is $T\kern-1pt$ normal whenever\/
$T^n\!$ is normal for some}\/ $n$\kern1pt?

\section{Multiples of Unitary Operators}

From now on all Banach and Hilbert spaces are complex.\ The aim of this
section is to prove Theorem 5.2, which gives an ultimate condition for a
normaloid operator to be a multiple of a unitary operator, and to disclose
some useful corollaries of it.\ To begin with, we encapsulate some necessary
basic results on invertibility.\

\vskip6pt\noi
{\bf Remark 5.1.}\
Recall that the inverse of an invertible Banach-space operator is an
operator~(i.e., it is bounded) by the Inverse Mapping Theorem.\

\vskip4pt\noi
{\rm(a)}
If $T\kern-1pt$ is an operator on a Banach space, then
${\sigma(T)^n=\sigma(T^n)}$ for every positive integer $n$, and so
\vskip4pt\noi
\centerline {\it $T^n\kern-1pt$ is invertible if and only if\/ $T\kern-1pt$ is
invertible.}\
\vskip4pt\noi
{\rm(b)}
Suppose $T\kern-1pt$ is invertible.\ Recall that (since ${I=T\kern1pt T^{-1}}$)
$$
1\le\|T\|\,\|T^{-1}\|.
$$
{\rm(c)}
So ${\|T^k\|^{-1/k}\!\le\|T^{-k}\|^{1/k}}$ for every
${k\kern-1pt\ge\kern-1pt1}.$ Thus, by the Gelfand--Beurling formula,
which says that ${r(T)=\lim_k\|T^k\|^{1/k}}$ for every operator $A$, we get
$$
r(T)^{-1}\le r(T^{-1})
\qquad\big(\hbox{equivalently}, \quad r(T^{-1})^{-1}\le r(T)\big).
$$
Normaloidness for both $T\kern-1pt$ and $T^{-1}\kern-1pt$ does not imply that
the above inequalities become identities (as we will see in Proposition 5.6(c)
below.)\                                                 \hfill$\blacksquare$

\vskip6pt
Let $T\kern-1pt$ be an operator on a Hilbert space.\

\vskip6pt\noi
{\bf Theorem 5.2.}\
{\it If\/ $T\kern-1pt$ is normaloid and\/ $T^n\kern-1pt$ is invertible with\/
${\|T^n\|^{-1}\!=\kern-1pt\|T^{-n}\|}$ for some positive integer\/ $n$, then\/
$T\kern-1pt$ is a multiple of a unitary operator}\/.\

\proof
Since $T\kern-1pt$ is invertible, ${\|T\|\ne0}.$
\vskip4pt\noi
(a)
Set ${S=\kern-1ptT/\|T\|}.$ As $T\kern-1pt$ is normaloid, $S$
is a normaloid contraction with ${\|S^k\|=1}$ for every positive integer $k$
(because ${S^k\kern-1pt=T^k/\|T\|^k\kern-1pt=T^k/\|T^k\|}\kern.5pt$).\
\vskip6pt\noi
(b)
Since $T^n$ is invertible, $S^n$ is invertible, and so is $S.$
\vskip6pt\noi
(c)
Since ${\|T^n\|^{-1}\!=\kern-1pt\|T^{-n}\|}$ and $T$ is normaloid, it 
follows that ${\|S^{-n}\|=1}.\;$ Indeed,
$$
\|S^{-n}\|=\|(T/\|T\|)^{-n}\|=\|T^{-n}\|/\|T\|^{-n}
=\|T^{-n}\|\,\|T\|^n=\|T^n\|^{-1}\,\|T^n\|=1.
$$
(d)
Thus, according to Remark 5.1(b), we get ${\|S^{-1}\|=1}.\;$ In fact,
$$
1=\|S\|^{-1}\le\|S^{-1}\|=\|S^{n-1}S^{-n}\|\le\|S^{n-1}\|\,\|S^{-n}\|
=\|S^{-n}\|=1.
$$
(e)
Since ${\|S\|=\|S^{-1}\|=1}$, $S$ is unitary.
\vskip6pt\noi
Actually, {\it on a Hilbert space, an invertible contraction whose inverse is
a contraction is unitary} (see, e.g.,
\cite[Proposition 5.73(a,d)]{EOT}$\kern.5pt$), which is a particular case of a
classical result that says that {\it an invertible power-bounded operator with
a power-bounded inverse is similar to a unitary operator} --- cf.\
\cite[Theorem 1]{Nag}$\kern.5pt$.\
\vskip6pt\noi
(f)
And so ${T=\|T\| S}$ is a multiple of a unitary operator.\               \qed

\vskip6pt
Immediate corollaries of Theorem 5.2 read as follows.\

\vskip6pt\noi
{\bf Corollary 5.3.}\
{\it If\/ $T\kern-.5pt\kern-1pt$ is normaloid and\/ $\kern-.5ptT^n\kern-1.5pt$
is a multiple of a unitary operator~for some positive integer\/ $\kern-.5ptn$,
then\/ $\kern-.5ptT\kern-1.5pt$ is a multiple of a unitary operator\/
$\kern-.5pt(\kern-.5pt$and thus~\hbox{normal}\/$)$}\/.\

\vskip6pt\noi
{\bf Corollary 5.4.}\
{\it If\/ $T\kern-1pt$ is normaloid and\/ $T^n\kern-1pt$ is a scalar operator
for some positive integer\/ $n$, then\/ $T\kern-1pt$ is a multiple of
a unitary operator\/ $($and thus normal\/$)$}\/.\

\vskip6pt
By a scalar operator we mean a multiple of the identity.\ The particular
case~in Corollary 5.4 has been considered in \cite[Lemma 4.1]{SK2}.\ Note that
if $T^n\kern-1pt$ is a zero multiple of a unitary or of the identity, the
above corollaries hold trivially with $T\kern-1pt$ being the null operator
because the only normaloid nilpotent operator is the null \hbox{operator}.\
Otherwise, $T^n\kern-1pt$ is invertible (and so is $T\kern-1pt$), as required
in Theorem 5.2.\ It is also worth noting that in case of Corollary 5.4, the
resulting operator~$T\kern-1pt$~is~not~necessari\-ly~a~multiple of the
identity (i.e., a scalar).\ For instance,
${T\kern-1pt=\big(\smallmatrix{0 & 1 \cr
                               1 & 0 \cr}\big)}$
and
${T^2\kern-1pt=\big(\smallmatrix{1 & 0 \cr
                                 0 & 1 \cr}\big)}.$

\vskip6pt
A version of the root problem asks for normaloid operators that become normal
when raised to an integer power.\ We consider below three properties of
normaloid operators that are required in the sequel, leading to two classes of
operators where normaloidness implies~\hbox{normality}.\

\vskip6pt\noi
{\bf Remark 5.5.}\
Let $T\kern-1pt$ be an operator on a Banach space.\
\vskip6pt\noi
{\rm(a)}
{\it If\/ $T\kern-1pt$ is normaloid, then\/ $T^k\kern-1pt$ is normaloid
for every integer}\/ ${k\kern-1pt\ge\kern-1pt1}.$
\vskip4pt\noi
(In fact,
${r(T)^k\kern-1.5pt=\kern-.5pt r(T^k)\kern-1pt\le\kern-1pt\|T^k\|
\kern-1pt\le\kern-1pt\|T\|^k\kern-1.5pt=\kern-.5ptr(T)^k\!}$,
so
${r(T^k)\kern-1pt=\kern-.5ptr(T)^k\kern-1.5pt=\kern-.5pt\|T\|^k
\kern-1.5pt=\kern-.5pt\|T^k\|}.$)
\vskip6pt\noi
{\rm(b)}
{\it On a $2$-dimensional space, normaloidness coincides with normality}\/.\ \\
(See, e.g., \cite[Theorem 3.3]{SK2}.)\
\vskip6pt\noi
{\rm(c)}
{\it For an involution, normaloidness coincides with normality}\/.\
\vskip4pt\noi
Theorem 5.2 can be used to supply a proof of this assertion.\ In fact, if
$T\kern-1pt$ is~an~involution (i.e., ${T^2\kern-1pt=I}$; equivalently,
${T^{-1}\kern-1pt=T}$), then ${r(T)^2=r(T^2)=r(I)=1}.$~Thus,
$$
1=r(T)\le\|T\|=\|T^{-1}\|.
$$
If, in addition, $T\kern-1pt$ is normaloid, then ${r(T)=\|T\|}$ so that
$$
\|T\|^{-1}=\|T^{-1}\|.
$$
Since $T$ is an invertible normaloid for which the above identity holds,
Theorem 5.2 with ${n=1}$ says that $T\kern-1pt$ is a multiple of a unitary
operator $U.$ As ${\|T\|=1}$, we get~$T=\gamma\,U$ with $\gamma$ in the unit
circle.\ $\kern-1pt$Then
${T^*\!=\overline\gamma\,U^*\!=\overline\gamma\,U^{-1}
\!=\overline\gamma\,\gamma^{-1}T^{-1}\!=T^{-1}\!=T}.$
$\kern-.5pt$So~$T\kern-1pt$ is a self-adjoint involution (i.e., a normal
involution).\ The converse is trivial$.$~Thus,~in addition to the equivalent
conditions in Section 2, $T\kern-1pt$ is a symmetry if and only if
\begin{description}
\item{$\kern3pt$\rm(vi)$\kern1pt$}
$\,T$ is a normaloid involution\/.                       \hfill$\blacksquare$
\end{description}
\vskip-2pt

\vskip6pt
Consider the statement in Theorem 5.2.\ Since $T\kern-1pt$ is an invertible
normaloid, the condition ${\|T^n\|^{-1}\!=\kern-1pt\|T^{-n}\|}$ for some $n$
is equivalent to ${\|T\|^{-n}\!=\kern-1pt\|T^{-n}\|}$ for such an $n.$ The
conclusion of Theorem 5.2 says that $T\kern-1pt$ is normal.\ If $T\kern-1pt$
is an invertible normal, then $T^{-1}\!$ is normal; so $T\kern-1pt$ and
$T^{-1}\!$ are normaloid, which implies by Remark 5.5(a) that
$\|T^k\|\kern-1pt=\kern-1pt\|T\|^k\!$ and
$\|T^{-k}\|\kern-1pt=\kern-1pt\|T^{-1}\|^k\!$ for every $k.$ In this case,
$$
\|T^n\|^{-1}\!=\kern-1pt\|T^{-n}\|
\kern-1pt\iff\kern-1pt
\|T\|^{-n}\!=\kern-1pt\|T^{-1}\|^n
\kern-2pt\kern-1pt
\iff
\|T\|^{-1}\!=\kern-1pt\|T^{-1}\|
\kern-2pt\iff\kern-1pt
1\kern-1pt=\kern-1pt\|T\|\kern1pt\|T^{-1}\|.
$$
Proposition 5.6 below completes such an argument.\

\vskip6pt\noi
{\bf Proposition 5.6.}\
{\it Let\/ $T\kern-1pt$ be an invertible operator on a Banach space
and consider the following assertions}\/.\
\begin{description}
\item{$\:$(i)$\:$}
${\|T\|^{-1}=\|T^{-1}\|}.$
\vskip4pt\noi
\item{$\kern.5pt$(ii)$\kern.5pt$}
 ${\|T^k\|^{-1}=\|T^{-k}\|}\;$ {\it for every positive integer}\/ $k.$
\vskip4pt\noi
\item{(iii)}
${\|T^n\|^{-1}=\|T^{-n}\|}\;$ {\it for some positive integer}\/ $n.$
\vskip4pt\noi
\item{(iv)}
${\|T^k\|^{-1}=\|T^{-k}\|}\;$ {\it for every ${k\ge n}$, for some}\/ $n.$
\vskip4pt\noi
\item{$\kern1.5pt$(v)$\kern1.5pt$}
${r(T)^{-1}=r(T^{-1})}.$
\end{description}
\vskip4pt\noi
{\rm(a)}
{\it If\/ $\kern-.5pt T\kern-1.5pt$ is normaloid, then {\rm(i)}$\kern-1pt$
implies {\rm(ii)},$\kern-1pt$ and {\rm(i)}$\kern-1pt$ also implies that\/
$\kern-.5pt T^{-1}\kern-2pt$ is \hbox{normaloid}}\/.\
\vskip4pt\noi
{\rm(b)}
{\it If\/ both\/ $T\kern-1pt$ and\/ $T^{-1}\kern-1.5pt$ are normaloid, then
assertions\/ {\rm(i)} to\/ {\rm(v)} are equivalent}\/.\
\vskip4pt\noi
{\rm(c)}
{\it $T\kern-1pt$ and\/ $T^{-1}\kern-1.5pt$ both being normaloid does not
imply any of the assertions\/ {\rm(i)} to}\/ {\rm(v)}.\

\proof
(a)
If (i) holds, then
${\|T^{-k}\kern-1pt\|\le\|T^{-1}\|^k\kern-2pt=\|T\|^{-k}}\kern-2pt$, and hence
${\|T\|^k\|T^{-k}\|\le1}$, for every ${k\ge1}.$ Thus, by Remark 5.1(b), if
$T\kern-1pt$ is normaloid, then
${1\le\|T^k\|\,\|T^{-k}\|}={\|T\|^k\|T^{-k}\|\le1}.$ Hence,
${\|T^k\|\,\|T^{-k}\|=1}$; that is, ${\|T^{-k}\|=\|T^k\|^{-1}}\!.$
So~(i)~$\Rightarrow$~(ii).\
\vskip4pt\noi
Also, if $T\kern-1pt$ is normaloid and (i) holds, then (ii) holds as we saw
above, and so (v) holds by the Gelfand--Beurling formula.\ Thus, since
${\|T\|^{-1}\!=r(T)^{-1}}\!$ (for $T\kern-1pt$~is~normaloid),
\vskip3pt\noi
$$
\|T^{-1}\|
\;{{{\buildrel_{\scriptstyle\rm(i)}\over=}}}\;
\|T\|^{-1}=r(T)^{-1}
\;{{{\buildrel_{\scriptstyle\rm(v)}\over=}}}\;
r(T^{-1}),
$$
\vskip4pt\noi
and so (i) implies that $T^{-1}\kern-1pt$ is normaloid whenever $T\kern-1pt$
is.\
\vskip4pt\noi
(b)
In general, for every invertible operator, (ii) $\Rightarrow$
$\big\{$(iii) and (iv)$\big\}$ trivially, and (iv) $\Rightarrow$ (v) by the
Gelfand--Beurling formula.\ If $T\kern-1pt$ is normaloid, then
(i) $\Rightarrow$ (ii) by item (a).\ Moreover, if $T\kern-1pt$ and
$T^{-1}\kern-1pt$ are normaloid, then (iii) $\Rightarrow$ (i).\ Indeed,
\vskip3pt\noi
$$
\|T\|^{-n}=\|T^n\|^{-1}
\;{{{\buildrel_{\scriptstyle\rm(iii)}\over=}}}\;
\|T^{-n}\|=\|(T^{-1})^n\|=\|T^{-1}\|^n,
$$
\vskip5pt\noi
so that ${\|T\|^{-1}=\|T^{-1}\|}.$ Finally,	if $T\kern-1pt$ and
$T^{-1}\kern-1pt$ are normaloid, then (v) $\Rightarrow$ (i).\

\vskip4pt\noi
(c)
For instance, if ${T\!=\kern-.5pt\frac{1}{2}I\oplus 2I}$, then each of (i)
up to (v) fails.                                                         \qed

\vskip6pt\noi
{\bf Corollary 5.7.}\
{\it Let\/ $T\kern-1pt$ be an invertible normaloid operator on a Hilbert
space.\ If\/ ${\|T\|^{-1}\kern-1pt=\|T^{-1}\|}$, then\/ $T\kern-1pt$ is a
multiple of a unitary operator, thus~normal}\/.\

\proof
$\!$Take any positive integer $n.$ If $T\kern-1pt$ is invertible, normaloid
and ${\|T\|^{-1}\!=\kern-1pt\|T^{-1}\|}$, then Proposition 5.6(a) ensures that
${\|T^n\|^{-1}\kern-1pt=\|T^{-n}\|}.$ So $T\kern-1pt$ is a multiple of a
unitary operator by Theorem 5.2.\                                        \qed

\vskip6pt
Corollary 5.7, or Proposition 5.6(a), ensures that $T^{-1}$ is normaloid as
well.\

\vskip6pt\noi
{\bf Corollary 5.8.}\
{\it Suppose\/ $T\kern-1pt$ is an invertible normaloid operator on a Hilbert
space with a normaloid inverse.\ The following assertions are equivalent}\/.
\begin{description}
\item{\rm(a)}
${r(T)^{-1}=r(T^{-1})}.$
\vskip4pt\noi
\item{\rm(b)}
${\|T\|^{-1}=\|T^{-1}\|}.$
\vskip4pt\noi
\item{\rm(c)}
{\it ${\|T^n\|^{-1}=\|T^{-n}\|}$ for some positive integer}\/ $n.$
\vskip4pt\noi
\item{\rm(d)}
{\it $T\kern-1pt$ is a multiple of a unitary operator}\/.
\end{description}

\proof
If $T$ and $T^{-1}$ are normaloid, then Proposition 5.6(b) ensures that (a) to
(c) are equivalent, (d) trivially implies (c), and (c) implies (d) by
Theorem 5.2.\                                                            \qed

\vskip6pt
Although being an invertible normaloid with a normaloid inverse does not
imply any of the above assertions, in this case the assertions are
equivalent:\ they all hold or they all fail together.\ Thus, if we add an
extra assumption to Corollary 5.7, namely, the inverse of $T\kern-1pt$ is
normaloid, then the central condition, ``${\|T\|^{-1}=\|T^{-1}\|}$'', can
be replaced by any of the equivalent conditions in Corollary 5.8.\

\section{The Root Problem for Completely Normaloid Operators}

Let $T\kern-1pt$ be a normaloid operator on a Hilbert space.\ Suppose all
parts of $T\kern-1pt$ are normaloid.\ If $T^n\!$ is normal for some $n$, then
$T\kern-1pt$ is normal.\ In other words:

\vskip6pt\noi
{\bf Theorem 6.1.}\
{\it On an arbitrary Hilbert space, if\/ an operator $T$ is
completely normaloid and\/ $T^n\!$ is normal for some positive integer $n$,
then\/ $T\kern-1pt$ is normal}\/.\

\vskip6pt
Theorem 6.1 gives an affirmative answer to Question 4.2.\ 

\vskip6pt
First we establish in Lemma 6.3 a version of the above result for a separable
\hbox{Hilbert} space.\ But before proving Lemma 6.3, we need a few results for
direct integrals on separable Hilbert spaces, which are summarised below and
reflect the minimum required in the proof of Lemma 6.3.\ For a comprehensive
treatment on direct integrals, including all those properties, see
\cite[Part II, Chapters 1 and 2]{Dix} (for a brief review focusing on normal
operators, see \cite[Section 4.1]{Geo}; see also~\cite{AC}$\kern.5pt).$

\vskip6pt\noi
(a)
A separable Hilbert space $\H$ can be represented as a direct integral (see,
e.g., \cite[Theorem 14.2.1]{KR}) of a family
${\{\H_\lambda\}\!=\!\{\H_\lambda\}_{\lambda\in\Lambda}}\!$ of Hilbert spaces
$\H_\lambda$, referred to as a field of fibre spaces, with respect to a
measure space ${(\Lambda,\mu)}$, where $\mu$ is a positive Borel measure on
the $\sigma$-algebra of Borel subsets of a locally compact complete metric
space $\Lambda.$ In the present case, $\Lambda$ will be a compact set in
$\CC$; thus, $\mu$ is finite.\ A direct integral representation of $\H$ is
usually denoted by
$$
\H=\!\int_\Lambda^{_\oplus}\H_\lambda d\mu.
$$
An element $x$ in $\kern-1pt{\int_\Lambda^{_\oplus}\H_\lambda d\mu}$ is a
function ${x\!:\Lambda\!\to\!\bigcup_\lambda\!\H_\lambda}$ such that
${x_\lambda\in\H_\lambda}$ for~each $\lambda$ in $\Lambda$, and
${\int_\Lambda^{_\oplus}\H_\lambda d\mu}$ is the collection of all these
functions ${x=\{x_\lambda\}=\{x_\lambda\}_{\lambda\in\Lambda}}$ for which
${\int_\Lambda\|x_\lambda\|^2d\mu<\infty}$, satisfying some additional
measurability conditions including the separability of each fibre space
$\H_\lambda$ (see, e.g.,
\cite[p$.\,$164, Definition 1, Section 3, Chapter$\;$1, Part$\;$II]{Dix}).\
Direct integrals extend the notion of direct sums.\

\vskip6pt\noi
(b)
A decomposable operator $T\kern-1pt$ on
${\H\kern-1pt=\!\int_\Lambda^{_\oplus}\H_\lambda d\mu}$ is one such that
${Tx_\lambda\kern-1pt=T_\lambda x_\lambda\kern-1pt}$ in $\H_\lambda$
($\mu$-a.e.)\ for a collection $\{T_\lambda\}$ of operators
${T_\lambda\!:\H_\lambda\!\to\kern-1pt\H_\lambda}$ with 
${\mu\hbox{-}\esup_{\lambda\in\Lambda}\|T_\lambda\|\kern-1pt<\!\infty}$, so
that ${T\kern-1ptx\kern-1pt=\kern-1pt\{T_\lambda x_\lambda\kern-.5pt\}}$ lies
in ${\!\int_\Lambda^{_\oplus}\kern-1pt\H_\lambda d\mu}$ for every
${x=\{x_\lambda\kern-.5pt\}}$ in
${\!\int_\Lambda^{_\oplus}\kern-1pt\H_\lambda d\mu}$ (satisfying additional
measurability conditions).\ A decomposable operator
${T\kern-1pt\!:\H\kern-1pt\to\kern-1pt\H}$ is denoted~by
$$
T\kern-1pt=\!\int_\Lambda^{_\oplus}T_\lambda d\mu,
$$
which is referred to as the direct integral of $T\kern-1pt.$ The induced
uniform norm of $T$ is
\vskip4pt\noi
(b$_1$)
${\|T\|=\mu\hbox{-}\esup_{\lambda\in\Lambda}\|T_\lambda\|}.$
\vskip4pt\noi
If $S$ on ${\H=\!\int_\Lambda^{_\oplus}\H_\lambda d\mu}$ is decomposable with
${S=\int_\Lambda^{_\oplus} S_\lambda d\mu}$, then
\vskip4pt\noi
(b$_2$)
${S=T\kern-1pt}$ if and only if ${S_\lambda=T_\lambda}$ $\mu$-a.e.,
\vskip4pt\noi
(b$_3$)
$TS$ is decomposable and
${TS=\!\int_\Lambda^{_\oplus} T_\lambda S_\lambda d\mu}$,
\vskip4pt\noi
(b$_4$)
$S^*$ is decomposable, and
${S^*\!=\!\int_\Lambda^{_\oplus} S_\lambda^* d\mu}.$
\vskip4pt\noi
From (b$_2$), (b$_3$) and (b$_4$) we may infer the following further
properties.\
\vskip4pt\noi
(b$_5$)
For every positive integer $n$, $T^n\kern-1pt$ is decomposable, and
${T^n\kern-1pt=\!\int_\Lambda^{_\oplus} T_\lambda^n d\mu}.$
\vskip4pt\noi
(b$_6$)
A decomposable $T\kern-1pt$ is normal if and only if the operators
$T_\lambda\kern-1pt$ are normal $\mu$-a.e.\
\vskip8pt\noi
(c)
If $N$ is a normal operator on $\H$, then it is a special kind of decomposable
operator, namely, a diagonalisable operator, where every $N_\lambda$ is a
scalar operator,
$$
N=\!\int_\Lambda^{_\oplus}\lambda\,I_\lambda d\mu,
$$
\vskip-1pt\noi
with ${\Lambda=\sigma(N)}$, the spectrum of the normal operator $N$, where
$\lambda$ lies in $\sigma(N)$ and each $I_\lambda$ is the identity operator
on $\H_\lambda.$ An important property of the von~\hbox{Neumann} algebra
generated by a normal operator $N$ on $\H$ is that every operator in its
commu\-tant is decomposable on
${\H=\int_{\sigma(N)}^{_\oplus}\H_\lambda\,d\mu}.$ In fact, an operator is
decomposable~if and only if it commutes with a diagonalisable operator (see,
e.g., \cite[p.$\,$188, Corollary, Section 5, Chapter 2, Part II]{Dix}).\
Hence,
\vskip4pt\noi
(c$_1$)
if $T\kern-1pt$ commutes with $N$, then
${T\kern-1pt=\!\int_{\sigma(N)}^{_\oplus}T_\lambda\,d\mu}$
with each $T_\lambda$ acting on $\H_\lambda.$

\vskip6pt
That having been said, assertions about properties of the operators
$T_\lambda$ of a de\-com\-posable operator
${T\kern-1pt=\!\int_\Lambda^{_\oplus} T_\lambda d\mu}$ are to be understood
$\mu$-almost everywhere.\ In particular, we say that the operators $T_\lambda$
are normaloid ($\mu$-a.e$.$) if ${\|T_\lambda^n\|=\|T_\lambda\|^n}$ for every
positive integer $n$ or, equivalently, ${\sigma(T_\lambda)=\|T_\lambda\|}$ or,
still equivalently, if ${w(T_\lambda)=\|T_\lambda\|}$, where these identities
are understood $\mu$-almost everywhere over $\Lambda$.\

\vskip6pt
The following auxiliary result is required in the proof of Lemma 6.3.\

\vskip6pt\noi
{\bf Proposition 6.2.}\
{\it Let\/ ${\H=\!\int_\Lambda^{_\oplus}\H_\lambda d\mu}$ be a direct integral
of a separable Hilbert space and let\/
${T\kern-1pt=\!\int_\Lambda^{_\oplus}T_\lambda d\mu}$ be a decomposable
operator on\/ $\H.$ If\/ $T\kern-1pt$ is completely normaloid, then\/
$T_\lambda$ is normaloid}\/ $\mu$-a.e.\

\proof
Recall that (i) an operator $A$ is normaloid if and only if ${\|A\|=w(A)}$,
where $w(A)$ stands for the numerical radius of $A$ (cf. Remark 3.1), (ii)
if ${A=\!\int_\Lambda^{_\oplus}T_\lambda d\mu}$ is a decomposable operator,
then ${\|A\|=\mu\hbox{-}\esup_{\lambda\in\Lambda}\|A_\lambda\|}$ (cf.\
(b$_1)\kern.5pt)$, and (iii) if $\Lambda$ is a locally compact Hausdorff space
and, in particular, if $\Lambda$ is a compact set in $\CC$ (as we have assumed
above), then ${w(A)=\esup_{\lambda\in\Lambda}\,w(A_\lambda)}$
\cite[Lemma]{Cho}.

\vskip6pt\hskip-4pt
Let $\Lambda'$ be the set of all ${\lambda\in\Lambda}$ for which $T_\lambda$
is not normaloid.\ That is,
$$
\Lambda'=\big\{\lambda\in\Lambda\!:\|T_\lambda\|>w(T_\lambda)\big\}.
$$
As the operator $T$ is decomposable, the maps ${\lambda\mapsto\|T_\lambda\|}$
and ${\lambda\mapsto w(T_\lambda)}$ are measur\-able, and so is the map
${\lambda\mapsto\|T_\lambda\|\kern-1pt-\kern-1pt w(T_\lambda)}$, which ensures
that $\Lambda'$ is a \hbox{measurable}~set.\

\vskip6pt\hskip-4pt
Suppose ${\mu(\Lambda')>0}.$ We can express $\Lambda'$ as the countable union
${\Lambda'=\bigcup_{k=1}^\infty\Lambda'_k}$,~where
$$
\Lambda'_k
=\big\{\lambda\in\Lambda\!:\|T_\lambda\|>w(T_\lambda)+\smallfrac{1}{k}\big\}
$$
It is clear that each $\Lambda'_k$ is also measurable and
${\Lambda'\subset\Lambda'_k}.$ Since ${\mu(\Lambda')>0}$, there exists an
integer $k$ such that ${\mu(\Lambda'_k)>0}.$ Let $T_k$ denote the restriction
of $T\kern-1pt$ to the subspace
${\H_k\kern-1pt=\!\int_{\Lambda'_k}^\oplus\H_\lambda d\mu}$ of $\H.$ As $\H_k$
reduces $T\kern-1pt$ (invariant to $T\kern-1pt$ and $T^*\!$ ---
cf$.$~(b$_4)\kern.5pt)$, and $T\kern-1pt$ is completely normaloid, the
operator $T_k$ is normaloid.\ Therefore,
$$
\|T_k\|=w(T_k).                                                   \eqno{(\S)}
$$
Moreover, as we also saw above
$$
\|T_k\|=\esup_{\lambda\in\Lambda'_k}\|T_\lambda\|
\quad \text{and}\quad
w(T_k)=\esup_{\lambda\in\Lambda'_k}w(T_\lambda).
$$
For every ${\lambda\in\Lambda'_k}$, we have the pointwise inequality
${\|T_\lambda\|>\omega(T_\lambda)+\frac{1}{k}}.$ Taking the essential
supremum over $\Lambda'_k$ yields
\begin{align*}
\|T_k\|
&=\esup_{\lambda\in\Lambda'_k}\|T_\lambda\|
\ge\esup_{\lambda\in\Lambda'_k}\big(w(T_\lambda)+\smallfrac{1}{k}\big)  \\
&=\big(\esup_{\lambda\in\Lambda'_k}\omega(T_\lambda\big)+\smallfrac{1}{k}
=\omega(T_k)+\smallfrac{1}{k}.
\end{align*}
Thus we get ${\|T_k\|\ge w(T_k)+\frac{1}{k}}$, which contradicts the identity
$(\S).$
\vskip6pt\hskip-4pt
Such a contradiction ensures that the assumption ${\mu(\Lambda')>0}$ is
false.\ Therefore, ${\mu(\Lambda')=0}$, which means that
${\|T_\lambda\|=w(T_\lambda)}$ for $\mu$-almost every ${\lambda\in\Lambda}.$
Thus, $T_\lambda$ is normaloid $\mu$-almost everywhere.\                 \qed

\vskip6pt\noi
{\bf Lemma 6.3.}\
{\it On a separable Hilbert space, if\/ an operator $T\kern-1pt$ is completely
normaloid and\/ $T^n\!$ is normal for some positive integer $n$, then\/
$T\kern-1pt$ is normal}\/.\

\proof
Let $T\kern-1pt$ be an operator on a separable Hilbert space $\H.$ Take an
arbitrary positive integer $n.$ If $T^n\kern-1pt$ is a normal operator on 
${\H=\!\int_{\sigma(T^n)}^{_\oplus}\!\H_\lambda d\mu}$, then, as we saw in
(c), it is a diagonalisable operator,
$$
T^n=\!\int_{\sigma(T^n)}^{_\oplus}\lambda\,I_\lambda d\mu.
$$
Since $T\kern-1pt$ commutes with the normal operator $T^n\!$, (c$_1$) says
that $T\kern-1pt$ is a decomposable operator on
${\H=\!\int_{\sigma(T^n)}^{_\oplus}\!\H_\lambda d\mu}$, which means that
$$
T\kern-1pt=\!\int_{\sigma(T^n)}^{_\oplus} T_\lambda d\mu,
$$
with each $T_\lambda$ acting on each $\H_\lambda.$ Hence, by (b$_5$),
$$
T^n=\!\int_{\sigma(T^n)}^{_\oplus}T_\lambda^nd\mu,
$$
\vskip-4pt\noi
Thus, according to (b$_2$), $\mu$-a.e.,
$$
T_\lambda^n=\lambda\,I_\lambda.
$$
\vskip2pt\noi
Now suppose $T\kern-1pt$ is completely normaloid so that $T_\lambda$ is
normaloid $\mu$-a.e.\ by Proposition 6.2.\ Then the above identity and
Corollary 5.4 ensure that $T_\lambda$ is normal $\mu$-a.e.\ Thus
${T\kern-1pt=\!\int_{\sigma(T^n)}^{_\oplus}T_\lambda d\mu}$ is normal on
$\H$ according to (b$_6$).\                                             \qed

\vskip6pt
To prove the general case where $\H$ is not necessarily separable, we borrow
an approach introduced in \cite[Proof of Lemma 3.1]{Pie}.\

\vskip6pt\noi
{\it Proof of Theorem 6.1}\/.\
Suppose $T\kern-1pt$ is a completely normaloid operator on a Hilbert space
$\H$ such that $T^n$ is normal for some positive integer $n.$ Take
an arbitrary ${x\in\H}.$ With $\NN$ denoting the set of all positive integers,
${\NN_0\kern-1pt=\NN\cup\kern-1pt\0}$, and with $\bigvee$ standing for
\hbox{closure}~of~span,~set
$$
\H_x=\bigvee\big\{{T^*}^{i_k}T^{j_k}\cdots{T^*}^{i_1}T^{j_1}x\!:\,
(i_1,\dots,i_k),(j_1,\dots,j_k)\in\NN_0^k,\, k\in\NN\big\}\;\;\sse\;\;\H,
$$
which is a subspace of $\H.$ Since $\H_x$ is spanned by a countable set, it
is separable, thus a separable Hilbert space.\ Also, as is readily verified,
$\H_x$ is invariant to both $T\kern-1pt$ and $T^*\kern-1pt$, and so it reduces
$T\kern-1pt.$ Set ${T_x=T|_{\H_x}}\!$ on $\H_x.$ Since $T\kern-1pt$ is
completely normaloid and $\H_x$ reduces $T\!$, it follows that $T_x$ is
completely normaloid as well.\ Since $T^n$ is normal and $\H_x$ reduces
$T\kern-1pt$ (and so it reduces $T^n$), we get
\begin{eqnarray*}
T_x^{*n}T_x^n
&\kern-6pt=\kern-6pt&
{T|_{\H_x}\!\!}^{*n}\,{T|_{\H_x}\!\!}^n
=T^{*n}|_{\H_x}T^n|_{\H_x}
=T^{*n}T^n|_{\H_x} \\
&\kern-6pt=\kern-6pt&
T^nT^{*n}|_{\H_x}
=T^n|_{\H_x}T^{*n}|_{\H_x}
={T|_{\H_x}\!\!}^n\,{T|_{\H_x}\!\!}^{*n}
=T_x^nT_x^{*n}.
\end{eqnarray*}
Then $T_x^n$ is normal.\ Since $\H_x$ is separable, $T_x$ on $\H_x$ is
completely normaloid, and $T_x^n$ is normal, it follows from Lemma 6.3 that
$T_x$ is normal.\ Thus, as $\H_x$ reduces $T\kern-1pt$,
$$
T\kern1ptT^*x
=T|_{\H_x}\!\!\,T|_{\H_x}^*x
=T_xT^*_xx
=T_x^*T_xx
=T|_{\H_x}^*T|_{\H_x}\!\!\,x
=T^*Tx.
$$
As $x$ was arbitrarily taken from $\H$, the above identity ensures that
${T\kern1pt T^*\!=T^*T\kern-1pt}$, and therefore $T$ is normal.\         \qed

\section{Some Consequences of Theorem 6.1}

To avoid trivialities, we agree that all spaces have dimensions greater than
$1.$~The subspaces $\N(T)$ and $\R(T)^-\!$ of $\H$ stand for the kernel and
closure of the range of an operator $T\kern-1pt$ on a Hilbert space $\H.$
First we pose a few elementary auxiliary results.\

\vskip6pt\noi
{\bf Remark 7.1.}\
Let $T\kern-1pt$ be an operator on a Hilbert space.\
\vskip4pt\noi
(i) If\/ ${T\kern-1pt\in\hbox{{\sc CN}}}$, then
\begin{description}
\item{\rm(a)}
{\it every part of\/ $T\kern-1pt$ lies in\/ {\sc CN}},
\end{description}
by the very definition of {\sc CN} (we have already applied this property),
and
\begin{description}
\item{\rm(b)}
{\it every power of\/ $T\kern-1pt$ lies in\/ {\sc CN}}.\
\end{description}
Indeed,
${T\kern-1pt\in\hbox{\sc CN}\limply\kern-1pt T^k\kern-1.5pt\in\hbox{\sc CN}}$
for every positive integer $k$ since powers of normaloid operators are
normaloid.\ The converse fails:\
${T^2\!\in\hbox{\sc CN}\notlimply\kern-.5pt T\kern-1pt\in\hbox{\sc CN}}.$
(In fact, if ${T\kern-1pt=L\oplus I}$ with $L$ being a nilpotent with
${\|L\|\kern-1pt>\kern-1.5pt1}$, then $T^2$ is normal, thus
${T^2\!\in\hbox{\sc CN}}$ trivially, but $T$ is not normaloid.)\
\vskip4pt\noi
(ii) If $T^k\kern-1pt$ is normal for some positive integer $k$, then the
subspaces
\begin{description}
\item{\rm(c)}
{\it $\R(T^k)^-\!$ and\/ $\N(T^k)$ reduce\/ every operator $F$ that commutes
with $T^k$},
\end{description}
by the Fuglede--Putnam Theorem:\ if an operator commutes with a normal
operator $N$, then it commutes with $N^*\!.$ (As $F$ commutes with
$T^k\kern-2pt$,
$\kern1pt F\kern1pt T^{*k}\kern-2pt=\kern-1pt T^{*k}F\kern-1pt$, so
$T^kF^*\kern-2pt=\kern-1pt F^*T^k\kern-2pt$, hence $\N(T^k)$ and
$\R(T^k)^-\!$ are $F^*\kern-1pt $-invariant, and, trivially
$F$-invariant.)                                \hfill$\blacksquare$

\vskip6pt\noi
{\bf Corollary 7.2.}\
{\it Let\/ $T\kern-1pt$ be an operator on a Hilbert space\/ $\H$, and take an
arbitrary integer\/ ${k\ge2}.$ Suppose\/ $T^k\!$ and all its parts are
normaloid\/ $($i.e., suppose\/ ${T^k\kern-1.5pt\in\hbox{\sc CN}}).$ If\/
$T^n\kern-1.5pt$ is normal for some positive integer\/ $n$,~then\/
$$
T=L\oplus S
\quad\hbox{on}\quad
\H=\N(T^n)\oplus\N(T^n)^\perp,
$$
where\/ $L$ is a nilpotent operator on\/ $\N(T^n)$ of index\/
${j\le\min\{k,n\}}$,\/ $S$ on\/ $\N(T^n)^\perp\!$ is such that\/
$S^k\kern-1.5pt$ is normal\/ $($as well as $S^n)$, and any of these
parts~may~be~absent}\/.\

\proof
Since $T^n$ is normal, and since $\N(T^n)$ reduces $T\kern-1pt$ (cf.\
Remark 7.1.(c)$\kern.5pt$),
$$
T=L\oplus S
\quad\hbox{on}\quad
\H=\N(T^n)\oplus\N(T^n)^\perp\!,
$$
where ${L=T|_{\N(T^n)}}$ and ${S=T|_{\N(T^n)^\perp\!}}$ are parts of
$T\kern-1pt.$ Thus, ${L^n\kern-1pt=(T|_{\N(T^n)})^n}={T^n|_{\N(T^n)}=O}$, so
that $L$ is nilpotent of index ${j\le n}$, and
${S^n\kern-1pt=(T|_{\N(T^n)^\perp})^n}={T^n|_{\N(T^n)^\perp}\!}$ is normal.\
Take an arbitrary ${k\ge1}$, and suppose ${T^k=L^k\oplus S^k}$ lies in
{\sc CN} so that $L^k$ and $S^k$ are again in {\sc CN} (cf.\
Remark 7.1(a)$\kern.5pt$).\
\vskip6pt\noi
(i)
Since $L^k$ is normaloid, and recalling that $L^k$ is nilpotent because $L$
is nilpotent, we get ${L^k=O}$, and so $L$ is nilpotent of index ${j\le k}.$
\vskip6pt\noi
(ii)
Since $S^k$ lies in {\sc CN} and since $(S^k)^n=(S^n)^k$ is normal because
$S^n$ is normal, Theorem 6.1 ensures that $S^k$ is normal.\              \qed

\vskip6pt
For ${k=1}$, Corollary 7.2 collapses to Theorem 6.1, where ${L=O}.$

\vskip6pt
There exist irreducible (thus not normal) involutions; sample:\
$T\kern-1.5pt={\big(\smallmatrix{1 & 1    \cr
                                 0 & \!-1 \cr}\big)}.$
But involutions, in spite of having a normal square trivially, are normaloid
if and only if they are normal.\ This and Theorem 6.1 suggest the following
question.\

\vskip6pt\noi
{\bf Question 7.3.}\
{\it Is there an irreducible normaloid operator with a normal square}\/\,?

\vskip6pt\noi
{\bf Answer 7.4.}\
Consider the following result from \cite[Theorem 1]{RR}.\
\vskip6pt\noi
{\narrower
{\it A Hilbert-space operator $T\kern-1pt$ is the square root of a normal
operator if and only if it is of the form\/
$T\kern-1pt={A\oplus\big(\smallmatrix{B & \kern4pt C  \cr
                                      O & \kern-1pt-B \cr}\big)}$,
where\/ $A$ and\/ $B$ are normal operators and\/ $C$ is positive\/ $($i.e.,
nonnegative and injective\/$)$ and commutes with\/ $B$.}
\vskip0pt}
\vskip6pt\noi
Suppose $T^2\kern-1pt$ is normal.\ So we can apply the above result.\ If
$T\kern-1pt$ is irreducible, then one of the above two parts is absent.\ If
part
$\big(\smallmatrix{B & \kern4pt C  \cr
                   O & \kern-1pt-B \cr}\big)$
is absent, then ${T\kern-1pt=A}$ is an \hbox{irreducible} normal operator,
which is a contradiction.\ $\kern-1pt$This implies that part~$A$~is~absent.\
$\kern-1pt$Thus,
${T\kern-1pt=\big(\smallmatrix{B & \kern4pt C  \cr
                               O & \kern-1pt-B \cr}\big)}$
is irreducible and
${T^2\!=\kern-1pt
\big(\smallmatrix{B^2 & O            \cr
                  O   & B^2\kern-1pt \cr}\big)=B^2\oplus B^2}\kern-1pt.$
Hence, $\|T^2\|=\|B^2\|=\|B\|^2$ because $B$ is normal.\ If $T\kern-1pt$ is
normaloid, then ${\|T^2\|=\|T\|^2}\kern-1pt$, so that $\|T\|=\|B\|$, which is
again a contradiction since ${C\kern-1pt\ne\kern-1pt O}$ ($C$ is
injective)$.$~\hbox{Therefore},
\goodbreak\vskip3pt\noi
$$
\hbox{\it if\/ $T^2\kern-1pt$ is normal and\/ $T\kern-1pt$ is irreducible,
then\/ $T\kern-1pt$ is not normaloid\/.}
$$
\vskip1pt

We conclude with two remarks along the lines of Answer 7.4 and Theorem 6.1.\

\vskip6pt\noi
{\bf Remark 7.5.}\
Let $T\kern-1pt$ be an operator on an arbitrary Hilbert space.\ If it is in
{\sc CN}, then it is normaloid, but it may be irreducible.\ Since it is
normaloid, all powers of it are normaloid, but if it is irreducible, none are
normal, as we will see in item (b) below.\ First, we verify a generalisation
of Answer 7.4.\
\vskip6pt\noi
(a)
{\it If\/ $T\kern-1pt$ is normaloid and\/ $T^n\kern-1.5pt$ is normal for some
integer\/ ${n\kern-1pt\ge\kern-1pt1}$, then\/ $T\kern-1pt$ is reducible}\/.
\vskip6pt\noi
Indeed, let $T\kern1pt$ be a normaloid operator, take an arbitrary integer
${n\ge1}$, and suppose $T^n\kern-1pt$ is normal.\ If $T^n\kern-1pt$ is scalar,
then Corollary 5.4 ensures that $T\kern-1pt$ is normal and so reducible.\
Thus, suppose the normal operator $T^n\kern-1pt$ is not scalar.\ Since
$T^n\kern-1.5pt$ commutes with~$T\kern-1pt$, the Fuglede--Putnam Theorem
(as applied in Remark 7.1(c)$\kern.5pt$) ensures that $T\kern-1pt$ commutes
with $T^{*n}\kern-1pt.$ So there exists a nonscalar operator (viz.,
$T^n\kern-1pt$) that commutes with $T\kern-1pt$ and with $T^*\kern-1pt$, and
therefore $T\kern-1pt$ is reducible (see, e.g.,
\cite[Proposition IV.5.3]{Con}$\kern.5pt$).\
\vskip6pt\noi
Since $T\kern-1pt$ is an irreducible operator in {\sc CN} if and only if\/
$T\kern-1pt$ is an irreducible \hbox{normaloid}, a straightforward application
of item (a) reads as follows (compare with Answer~7.4).\
\vskip6pt\noi
(b)
{\it If\/ $T\kern-1pt$ is an irreducible operator in\/ {\sc CN}, then\/
$T^n\kern-1.5pt$ is not normal for all\/ ${n\ge1}$}\/,
\vskip6pt\noi
which also comes as a consequence of Theorem 6.1.\       \hfill$\blacksquare$

\vskip6pt\noi
{\bf Remark 7.6.}\
Let $T\kern-1pt$ be again a normaloid operator acting on a separable Hilbert
space $\H$, as in Lemma 6.3.\ If $T^n$ is normal, then $T\kern-1pt$ can be
decomposed into a countable orthogonal direct sum
$$
T={\bigoplus}_{k=0}^\infty T_k
\quad\;\hbox{on}\;\quad
\H={\bigoplus}_{k=0}^\infty\H_k,
$$
where ${T_0=T|_{\H_0}}$ is nilpotent and the remaining ${T_k=T|_{\H_k}}$ are
similar to normal operators $N_k$ on $\H_k$ (cf.\
\cite[Theorem 3.1 and the paragraph above it, p.$\,$139]{Gil}$\kern.5pt$).\
Since $T_0$ is normaloid, it is the null operator.\ Let $G$ be an invertible
(with a bounded inverse) operator on $\H_k$ for which ${T_k=G^{-1}N_kG}.$
Since ${\sigma(T_k)=\sigma(N_k)}$ for the normal operator $N_k$, there is a
scalar $\lambda$ such that ${T_k-\lambda I}$ and ${N_k-\lambda I}$
are invertible,~and
$$
T_{a,k}=\frac{(T_k-\lambda I)}{\|T_k-\lambda I\|}
=\frac{G^{-1}(N_k-\lambda I)G}{\|T_k-\lambda I\|}
=G^{-1}N_{a,k}G,
$$
is an invertible contraction similar to an invertible normal $N_{a,k}.$ The
fact that $T_{a,k}$~is an invertible contraction implies that
$$
\hbox{$T_{a,k}$ is similar to a unitary operator},
$$
say $T_{a,k}={W^{-1}UW}.\!$ Then ${N_{a,k}G\,W^{-1}=G\,W^{-1}U}$, which
implies by the Fuglede--Putnam Theorem that $N_{a,k}$ is unitary and the
operator $G\,T_{a,k}G^{-1}$ and its inverse are both power bounded.\ Hence
$T_{a,k}$ is a generalised scalar (cf.\
\cite[Theorem1.5.13, p.$\,$69]{LN}$\kern.5pt$). Evidently,
${\sigma(T_{a,k})\sse\TT}$, the unit circle.\            \hfill$\blacksquare$

\bibliographystyle{amsplain}

\end{document}